\documentclass{amsart}
\usepackage{epsfig}
\usepackage{amssymb}

\begin{document}

\newtheorem{thm}{Theorem}[section]
\newtheorem{lem}[thm]{Lemma}
\newtheorem{cor}[thm]{Corollary}

\theoremstyle{definition}
\newtheorem{defn}{Definition}[section]

\theoremstyle{remark}
\newtheorem{rmk}{Remark}[section]
\newtheorem{exa}{Example}[section]

\def\square{\hfill${\vcenter{\vbox{\hrule height.4pt \hbox{\vrule width.4pt
height7pt \kern7pt \vrule width.4pt} \hrule height.4pt}}}$}

\def\R{\mathbb R}
\def\Z{\mathbb Z}
\def\C{\mathbb C}
\def\H{\mathbb H}
\def\CP{\mathbb{CP}}

\newenvironment{pf}{\noindent {\it Proof:}\quad}{\square \vskip 12pt}

\title{Napoleon in isolation}

\author{Danny Calegari}
\address{Department of Mathematics \\ UC Berkeley \\ Berkeley, CA 94720}
\email{dannyc@math.berkeley.edu}

\maketitle

\begin{abstract}
Napoleon's theorem in elementary geometry describes how certain linear
operations on plane polygons of arbitrary shape always produce regular 
polygons. More generally,
certain triangulations of a polygon that tiles
$\R^2$ admit deformations which keep fixed the symmetry group of
the tiling. This gives rise to 
{\em isolation phenomena} in cusped hyperbolic $3$-manifolds, where 
hyperbolic Dehn surgeries on some collection of cusps leaves the 
geometric structure at some other collection of cusps unchanged.
\end{abstract}

\section{Geometric isolation}
\subsection{Definition}

Let $M$ be a complete, finite volume hyperbolic $3$-manifold with $n$
torus cusps, which we denote $c_1, \dots c_n$. The following definition is
found in \cite{wNaR93}:

\begin{defn}
A collection of cusps $c_{j_1},\dots c_{j_m}$ is {\em geometrically
isolated} from a collection $c_{i_1},\dots c_{i_n}$ if any hyperbolic
Dehn surgery on any collection of the $c_{i_k}$ leaves the geometric 
structure on all the $c_{j_l}$ invariant.
\end{defn}
Note that this definition is {\em not} symmetric in the collections
$c_{i_k}$ and $c_{j_l}$, and in fact there are examples which show that
a symmetrized definition is strictly stronger (see \cite{wNaR93}).

More generally, we can ask for some prescribed set of fillings on the
$c_{i_k}$ which leave the $c_{j_l}$ invariant. Generalized (non-integral)
hyperbolic surgeries on a cusp are a holomorphic parameter for the
space of all (not necessarily complete) hyperbolic structures on $M$
with a particular kind of allowable singularities (i.e. generalized
cone structures) in a neighborhood of the complete structure. Moreover, the
complex dimension of the space of geometric shapes on a complete cusp
is $1$. Consequently, for dimension reasons whenever $n>m$ there will be
families of generalized surgeries leaving the geometric structures on the
$c_{j_l}$ invariant. There is no particular reason to expect, however, that
any of these points will correspond to an {\em integral} surgery on the
$c_{i_k}$. When there is a $1$ complex dimensional holomorphic family
of isolated generalized surgeries which contains infinitely many integral 
surgeries, we say that we have an example of an {\em isolation phenomenon}.

Neumann and Reid describe other qualities of isolation in \cite{wNaR93}
including the following:
\begin{defn}
A collection of cusps $c_{j_1} \dots c_{j_m}$ is {\em strongly isolated}
from a collection $c_{i_1} \dots c_{i_m}$ if after any hyperbolic Dehn 
surgeries on any collection of the $c_{j_l}$, a further surgery on any
collection of the $c_{i_k}$ leaves the geometry of the (possibly filled)
cusps $c_{j_l}$ invariant.

A collection of cusps $c_{j_1} \dots c_{j_m}$ is {\em first-order isolated}
from a collection $c_{i_1} \dots c_{i_m}$ if the derivative of the
deformation map from generalized surgeries on the $c_{i_k}$ to the
space of cusp shapes on the $c_{j_l}$ vanishes at the complete structure.
\end{defn}

By using the structure of the $\Phi$ function defined in \cite{wNdZ85},
Neumann and Reid show that strong isolation and first order isolation
{\em are} symmetric relations. First order isolation can be restated
in terms of group cohomology, and is studied in some papers of Kapovich,
notably \cite{mK92}. 

In this paper we produce new constructions of isolation phenomena of
various qualities, both by extending or modifying known constructions,
and by introducing a conceptually original construction based on 
Napoleon's theorem in plane geometry. The sheer wealth of examples
that these techniques can produce, used in combination, strongly suggests
that instances of isolation phenomena are not isolated phenomena.

\subsection{Holomorphic rigidity}

Suppose we want to show that a cusp $d$ is isolated from a cusp $c$. Then
since the shape of the complete cusp $d$ depends holomorphically on the
generalized surgery on $c$, it suffices to show that infinitely many
surgeries on $c$ keep the structure of $d$ fixed, since these can only
accumulate at the complete (unfilled) structure, where we know the
function relating the structure on $d$ to the surgery on $c$ is regular.

We describe this holomorphic structure in more detail. Choose a 
meridian $m$ of a cusp $c$. Let $M_{(p,q)}$ be obtained by doing $(p,q)$
surgery on $c$. The hyperbolic structure on $M_{(p,q)}$ determines
a representation $\rho_{p,q}:\pi_1(M) \to PSL(2,\C)$. 
The image of $[m]$ under this representation has a well--defined {\em complex
length} $u$ in $\C/2\pi i \Z$, which is the logarithm of the ratio of the
eigenvalues of $\rho_{p,q}(m)$. We may choose a branch of the logarithm
so that the value $0$ corresponds to the complete structure. Then small
deformations of the hyperbolic structure on $M$ corresponding to generalized
surgeries on $c$ are parameterized in a
$2:1$ way by $u$. The set of Euclidean structures on a complete cusp
$d$ are parameterized by the complex orbifold 
$\mathcal{M}_1 = \H^2 / PSL(2,\Z)$, and
it is well--known that the map from $u \to \mathcal{M}_1$ 
is analytic, and regular at $0$. For more detail, see \cite{wNdZ85}.

\subsection{Rigid orbifolds}

A systematic study of isolation was initiated in \cite{wNaR93}.
Most of the examples constructed in \cite{wNaR93} of pairs of cusps $c_1,c_2$
which are geometrically isolated have the property that there is a totally
geodesic rigid triangle orbifold separating the two cusps. Such a separating
surface splits the manifold $M$ into two pieces $M_1,M_2$ where $c_i$ sits
as a cusp in $M_i$. Then a surgery on $c_i$ deforms only the piece $M_i$,
keeping the geometry of the splitting orbifold unchanged, and $M_i$ can then
be glued to $M_{i+1}$ to produce a complete hyperbolic structure on $M$.
One sees that the geometry of the entire piece containing the unfilled cusp
is unchanged by this operation, and therefore that each of the two cusps
is isolated from the other cusp.

If $M$ covers an orbifold $N$ containing a totally geodesic triangle orbifold
which separates cusps $c$ and $d$, then in $M$ any surgeries on lifts of 
$c$ which descend to $N$ will leave unaffected the structures on lifts of $d$.

\subsection{Rigid cusps}

A refinement of the construction above comes when the rigid orbifold is
boundary parallel. The square torus $\R^2/\Z\oplus i\Z$ and the hexagonal
torus $\R^2/\Z\oplus \frac {1+\sqrt{3}i} 2 \Z$ have rotational symmetries
of order $4$ and $6$ respectively. If these symmetries about a cusp $c$
can be made to extend over the entire manifold $M$, then any surgery preserving
these symmetries will keep the cusp $c$ square or hexagonal respectively.
Note that isolation phenomena produced by this method tend to be one--way.

\begin{exa}
Let $M=T^2 \times I - K$ where $K$ is a knot which has $4$--fold rotational
symmetry, as seen from either of the $T^2$ ends of $T^2 \times I$. Then
$M$ has $4$--fold rotational symmetry which is preserved after $(p,q)$
filling on the knot $K$. This symmetry keeps the shape of the two
ends of $T^2 \times I$ square after every surgery on $K$. In general,
surgery on either of the ends of the $T^2$ will disrupt the symmetry and
not exhibit any isolation.
\end{exa}

\subsection{Mutation}

An incompressible surface in a hyperbolic manifold does not need to be
rigid for certain topological symmetries to be realized as geometric
symmetries. In a finite volume hyperbolic $3$--manifold, an 
incompressible, $\partial$--incompressible 
surface $S$ without accidental parabolics
which is not the fiber of a fibration over $S^1$ is quasifuchsian, and
therefore corresponds to a unique point in $Teich(S) \times Teich(\bar{S})$.
For surfaces $S$ of low genus, the tautological curve over the
Teichm\" uller space of $S$ has certain symmetries which restrict to a
symmetry of each fiber --- that is, of each Riemann surface topologically
equivalent to $S$. There is a corresponding symmetry of the universal
curve over $Teich(\bar{S})$ and therefore of a quasifuchsian representation
of $S$. Geometrically, this means that one can cut along a minimal surface
representing the class of $S$ and reglue it by an automorphism which
preserves the geometric structure to give a complete 
nonsingular hyperbolic structure on a new manifold. Actually, one does
not need to know there is an equivariant minimal surface along which one
can cut --- one can perform the cutting and gluing at the level of limit
sets by using Maskit's combinations theorems (\cite{bM87}). This operation is
called {\em mutation} and is studied extensively by Ruberman, for instance
in \cite{dR87}.

If $S$ is a sphere with $4$ punctures, we can think of $S$ as a Riemann
sphere with $4$ points deleted. After a M\" obius transformation, we can
assume $3$ of these points are at $0,1,\infty$ and the $4$th is at the
complex number $z$. The full symmetry group $\mathsf{S}_4$ does not act
holomorphically on $S$ except in very special cases, but the subgroup 
consisting of $(12)(34)$ and its conjugates {\em does} 
act holomorphically. This group is known as the {\em Klein $4$--group}.

Geometrically, where we find an appropriate
$4$--punctured sphere in a hyperbolic manifold $M$, we can cut along the
sphere and reglue after permuting the punctures. This operation leaves
invariant the geometric structure on those pieces of the manifold that
do not meet the sphere. If $S$ is a spherical orbifold or cone manifold
with $4$ equivalent cone points, we can similarly cut and reglue.
The sphere along which the mutation is performed is called a
{\em Conway sphere}. The observation that mutation can be performed for
spherical orbifolds is standard: one can always find a finite manifold
cover of any hyperbolic orbifold, by Selberg's lemma. The spherical
orbifold lifts to an incompressible surface in this cover, and one can
cut and reglue equivariantly in the cover. In fact, local rigidity for
small cone angles developed by Hodgson and Kerckhoff in \cite{cHsK98} suggests
that this can be done just as easily for the cone manifold case, but this
is superfluous for our applications.

\begin{exa}
Let $M$ be a manifold with a single cusp $c$ which admits a $\Z/3\Z$ symmetry
that acts on $c$ as a rotation. This symmetry forces $c$ to be hexagonal.
Let $p$ be a point in $M$ not fixed by $\Z/3\Z$, and let $M'$ be obtained
by equivariantly removing $3$ balls centered at the orbit of $p$. Glue another
manifold with $\Z/3\Z$ symmetry to $M'$ along its spherical boundaries 
$S_1,S_2,S_3$ so that the symmetries on either side are compatible,
to make $M''$. Let $K$ be a $\Z/3\Z$-invariant knot in $M''$ which intersects
each of the spheres $S_i$ in $4$ points and which is sufficiently complicated
that its complement in $M''$ is atoroidal and the $4$--punctured spheres
$S_i$ are incompressible and $\partial$--incompressible. Perform a
mutation on $S_1$ which destroys the $\Z/3\Z$--symmetry to get a manifold
$N$ with two cusps which we refer to as $c$ and $d$, where $d$ corresponds
to the core of $K$. Then $c$ is hexagonal, since mutation does not
affect the geometric structure away from the splitting surface. Moreover,
for large integers $r$, 
any generalized $(r,0)$ surgery on $d$ will preserve the
fact that the $S_i$ are incompressible, $\partial$--incompressible
spheres with $4$ cone points of order $r$, and therefore we can undo the
mutation on $S_1$ to see that these surgeries do not affect the hexagonal
structure on $c$. But if a real $1$--parameter family of surgeries on
$d$ do not affect the structure on $c$, then {\em every} surgery on $d$
keeps the structure on $c$ fixed, so $c$ is isolated from $d$. 
Note that for general $(p,q)$ surgeries
on $d$, the spheres $S_i$ will be destroyed and the surgered manifold
will not be mutation--equivalent to a $\Z/3\Z$ symmetric manifold. Moreover,
for sufficiently generic $K$, there will be no rigid triangle orbifolds
in $N$ separating $c$ from $d$. As far as we know, this is the first
example of geometric isolation to be constructed that is not forced by
a rigid separating or boundary parallel surface.
\end{exa}   

We note in passing that mutation followed by surgery often preserves other
analytic invariants of hyperbolic cusped manifolds, such as volume
and (up to a constant), Chern--Simons invariants. Ruberman observed
in \cite{dR87} that if $K$ and a Conway sphere $S$ are {\em unlinked},
then a mutation which preserves this unlinking can actually be
achieved by mutation along the genus $2$ surface obtained by tubing
together the sphere $S$ with a tubular neighborhood of $K$. This
mutation corresponds to the hyperelliptic involution of a genus
$2$ surface. Since this surface is present and incompressible for
all but finitely many surgeries on $K$ and its mutant, the surgered
mutants are mutants themselves. If $S$ and $K$ are not unlinked, no
such hyperelliptic surface can be found, and in fact in this case
there is no clear relationship between invariants of the manifold obtained
after surgeries on the original
knot and on the mutant. Perhaps this makes the persistence of isolation
under mutation more interesting, since it shows that some of the effects of
mutation are persistent under surgeries where other effects are destroyed.

\section{Napoleon's theorem}

\subsection{Triangulations and hyperbolic surgery}

A finite volume complete but not compact hyperbolic $3$--manifold $M$
can be decomposed into a finite union of (possibly degenerate)
ideal tetrahedra glued together along their faces. An ideal tetrahedron
determines and is determined by its fourtuple of endpoints on the
sphere at infinity of $\H^3$. Identifying this sphere with $\CP^1$, 
we can think of the ideal tetrahedron as a $4$--tuple of complex
numbers. The isometry type of the tetrahedron is determined by the
{\em cross--ratio} of these $4$ points; equivalently, if we move three
of the points to $0,1,\infty$ by a hyperbolic isometry, the isometry type
of the tetrahedron is determined by the location of the $4$th point,
that is by a value in $\C -\lbrace 0,1 \rbrace$. This value is referred
to as the {\em simplex parameter} of the ideal tetrahedron. A combinatorial
complex $\Sigma$ 
of simplices glued together can be realized geometrically
as an ideal triangulation of a finite volume hyperbolic manifold (after 
removing the vertices of $\Sigma$) if
certain equations in the simplex parameters are satisfied. These equations
can be given explicitly by examining the links of vertices of the
triangulation. 

For a finite volume non--compact hyperbolic manifold, all the links of
vertices of $\Sigma$ should be tori $T_j$. There are induced triangulations
$\tau_j$ of these tori by the small triangles obtained by cutting off the tips
of the tetrahedra in $\Sigma$. Let $\Sigma = \bigcup_i \Delta_i$ and let
$z_i$ be a hypothetical assignment of simplex parameters to the tetrahedra
$\Delta_i$. Then the horoball sections of a hypothetical hyperbolic
structure on $\Sigma' = \Sigma - \text{vertices}$ are Euclidean tori
triangulated by Euclidean triangles, which are the horoball sections of
the $\Delta_i$. If the simplex parameter of an ideal tetrahedron is
$z$, the Euclidean triangles obtained as sections of the horoballs centered
at its vertices have the similarity type of the triangle in $\C$ with
vertices at $0,1,z$. A path in the dual $1$--skeleton to $\tau_j$ determines
a developing map of $\tau_j$ in $\C$: given a choice of the initial triangle,
there is a unique choice of each subsequent triangle within its fixed
Euclidean similarity type such that the combinatorial edge it shares
with its predecessor is made to geometrically agree with it.
The {\em holonomy} of a closed path in this dual skeleton is the Euclidean
similarity taking the initial triangle to the final triangle.

There are two necessary conditions to be met in order for
the hyperbolic structures on the $\Delta_i$ to glue up to give a hyperbolic
structure on $\Sigma'$. These conditions are actually sufficient in a
small neighborhood of the complete structure as described in \cite{wNdZ85}
and made rigorous in \cite{cPjP00}.

\begin{itemize}
\item{The {\em edge equations}: the holonomy around a vertex of
$\tau_j$ should be trivial.}
\item{The {\em cusp equations}: the holonomy around the meridian and
longitudes of the $T_j$ should be {\em translations}.}
\end{itemize}

These ``equations'' can be restated as identities of the form
$$\sum_i c_{ij} \ln(z_i) + d_{ij} \ln(1-z_i) = \pi i e_j$$
for some collection of integers $c_{ij},d_{ij},e_j$ 
and some appropriate choices
of branches of the logarithms of the $z_i$. 

A $(p,q)$ hyperbolic dehn surgery translates in this context to replacing the
cusp equations for some cusp by the condition that the holonomy around the
meridian and longitude give Euclidean similarities $h_m,h_l$ such that
$h_m^ph_l^q = \text{id}$. The analytic co--ordinates on Dehn surgery space
determined by the analytic parameters $z_i$ are holomorphically related to
the $u$ co--ordinates alluded to earlier. The geometry of a complete cusp
is determined by the translations corresponding to the holonomy of the
meridian and longitude.

This is all described in great detail in \cite{wNdZ85} and \cite{cPjP00}.

\subsection{Tessellations with forced symmetry}

In this section, we produce examples of isolation phenomena which do
not come from rigidity or mutation, but rather from the following theorem,
known as ``Napoleon's Theorem'' (\cite{hCsG67}):

\begin{thm}[``Napoleon's Theorem'']
Let $T$ be a triangle in $\R^2$. Let $E_1,E_2,E_3$ be three equilateral
triangles constructed on the sides of $T$. Then the centers of the $E_i$
form an equilateral triangle.
\end{thm}
\begin{pf}
Let the vertices of $T$ be $0,1,z \in \C$. The three centers of the $E_i$
are of the form $a_iz + b_i$ for certain complex numbers $a_i,b_i$
so the shape of the resulting triangle is
a holomorphic function of $z$. Let $\omega = (1+\frac i {\sqrt{3}})/2$. 
If $z$ is real and between $0$ and $1$ then
the center of $E_1$ splits the line between $0$ and $\omega$ in the
ratio $z,1-z$, the center of $E_2$ splits the line between $\omega$ and $1$
in the ratio $z,1-z$, and the center of $E_3$ is at $1-\omega$. But a
clockwise rotation through $\pi/3$ about $1-\omega$ takes the line between
$0$ and $\omega$ to the line between $\omega$ and $1$; i.e. it takes the
center of $E_1$ to the center of $E_2$. Thus the theorem is true for real
$z$ and by holomorphicity, it is true for all $z$. 
\end{pf}

\begin{figure}[ht]
\scalebox{.4}{\includegraphics{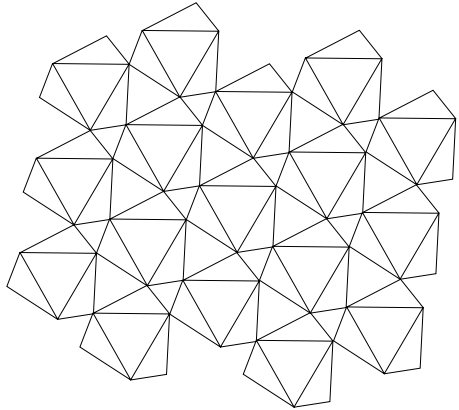}} 
\caption{Three isometric triangles and three equilateral triangles make up
a hexagon which tiles the plane with symmetry group $\mathsf{S}(3,3,3)$}
\end{figure}

Napoleon's theorem gives rise to an interesting phenomenon in plane
geometry: fix a triangle $T$. Then three triangles isometric to $T$
and three equilateral triangles with side length equal to the sides of
$T$ can be glued together to make a hexagon which tiles the plane with
symmetry group $\mathsf{S}(3,3,3)$ as in figure 1. The edge lengths of
the three equilateral triangles are equal to the three edge lengths of
$T$, and around each vertex the angles are the three interior angles of
$T$ together with three angles equal to $\frac {2\pi} 6$. It follows that
the hexagon in question exists and tiles the plane; to see that it has
the purported symmetry group, observe that the combinatorics of the
triangulation have a $3$--fold symmetry about the centers of the three
equilateral triangles bounding some fixed triangle of type $T$. These
three centers are the vertices of an equilateral triangle, by Napoleon's
theorem; it follows that the symmetry group of the tiling contains the
group generated by three rotations of order $3$ with centers at the
vertices of an equilateral triangle --- that is to say, the symmetry
group $\mathsf{S}(3,3,3)$. 

If we imagine that these triangles
are the asymptotic horoball sections of ideal tetrahedra going out a
cusp of a hyperbolic manifold, we see that appropriate deformations of
the tetrahedral parameters change the triangulation {\em but not the
geometry} of the cusp. For $T$ a horoball section of an ideal tetrahedron
with simplex parameter $z$, the holonomy around vertices is just
$$z \cdot \omega \cdot \frac {z-1} z \cdot \omega \cdot \frac 1 {1-z} \cdot
\omega = 1$$ 
where $\omega = \frac {1 + i\sqrt{3}} 2$ is the similarity type of an
equilateral triangle, and the holonomy around the meridian and longitude
are both translations by complex numbers $z_1,z_2$ whose ratio is
$\frac {1 + i\sqrt{3}} 2$.

\begin{figure}[ht]
\scalebox{.4}{\includegraphics{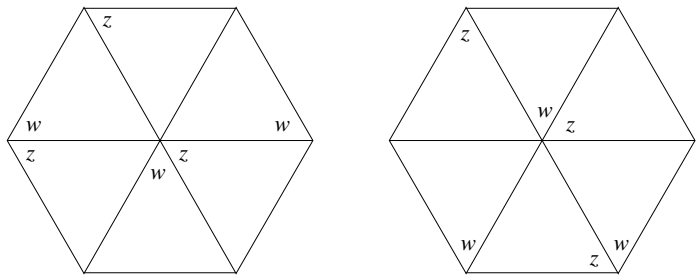}}
\caption{Two cusps with triangle parameters related as indicated in the
diagram can be deformed independently through (incomplete) affine structures.
We follow the conventions of \cite{wNdZ85}, so that if the vertices of
$\Delta$ are at $\lbrace 0,1,\infty,z \rbrace$, 
the edge marked with $z$ runs from $0$ to $\infty$ or from $1$ to $z$.}
\end{figure}

\begin{exa}
A scheme to parlay this theorem into isolation phenomena is given by
the following setup: we have two hexagons $H_1,H_2$ each tiled by six
equilateral triangles. We divide these twelve triangles into two
similarity types, corresponding to the triangles in $\C$ with vertices at
$\lbrace 0,1,z \rbrace$ and $\lbrace 0,1,w \rbrace$, 
and we choose a preferred vertex for each triangle corresponding to $0$, 
in the manner indicated in figure 2. For an arbitrary
choice of complex numbers $z$ and $w$,
the hexagons can be realized geometrically to give
affine structures on the tori obtained by gluing opposite sides of $H_1$
together and similarly for $H_2$. Initially set 
$z = w = \frac {1 + i\sqrt{3}} 2$. Deforming $z$ but keeping $w$ fixed
changes the affine structure on the first torus but leaves the structure
on the second torus unchanged. For, the combinatorial triangulation of
the second torus is exactly the triangulation of a fundamental domain of
the tessellation in figure 1. It follows that for $w = \frac {1 + i\sqrt{3}} 2$
and $z$ arbitrary, the universal cover of the second torus is tiled
by a tessellation with symmetry group $\mathsf{S}(3,3,3)$ and the torus is
therefore hexagonal.

Similarly, deforming $w$ and keeping $z = \frac {1 + i\sqrt{3}} 2$
changes the affine structure on the second torus but leaves the
structure on the first torus unchanged, since now we can identify the
combinatorial triangulation of the {\em first} torus with the triangulation
of a fundamental domain of the tessellation in figure 1.

This configuration of ideal tetrahedra with horoball sections equal to
the two cusps in this figure can be realized geometrically
by arranging six regular ideal tetrahedra in the upper--half space in a
hexagonal pattern with the common edge of the tetrahedra going from $0$ to
$\infty$. The pattern seen from infinity looking down is $H_1$, and the
pattern seen from $0$ looking up is $H_2$. The pictures are aligned so
that the real line has its usual orientation. Glue the twelve free faces of
the tetrahedra in such a way as to make $H_1$ and $H_2$ torus cusps. This
gives an orbifold $N$ whose underlying manifold is $T^2 \times I$,
with orbifold locus three arcs of cone angle
$2\pi/3$ each running between two $(3,3,3)$ triangle cusps arranged in the
obvious symmetrical manner.

Unfortunately, under surgeries of the cusps of $N$, the simplices do not
deform in the manner required by Napoleon's theorem. However, if we pass
to a $3$-fold cover, this problem can be corrected.

\begin{figure}[ht]
\scalebox{.3}{\includegraphics{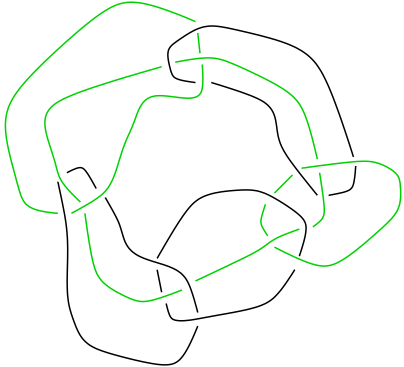}}
\caption{``Napoleon's $3$-manifold'' has a decomposition into ideal
tetrahedra such that surgery on a ``dark'' cusp deforms the triangulation
of the other ``dark'' cusps but keeps the shape of the cusp constant, as
guaranteed by Napoleon's theorem. A similar relation holds for the ``light''
cusps.}
\end{figure}

Let $L$ be the link depicted in figure $3$, and let $M = S^3 - L$. Then
$M$ admits a complete hyperbolic structure which can be decomposed into
$18$ regular ideal tetrahedra. It follows that $M$ is commensurable with
the figure $8$ knot complement. In fact, $M$ is the $3$-fold cover of $N$
promised above. 

Geometrically, arrange $18$ regular ideal tetrahedra
in the upper half-space with a common vertex at infinity so that a
horoball section intersects the collection the pattern depicted in figure
$4$. If we glue the $12$ external vertical sides of this collection of 
tetrahedra in the indicated manner, it gives a horoball section of one of the
cusps of $M$. It remains to glue up the $18$ faces of the ideal tetrahedra
with all vertices on $\C$. Figure $4$ has an obvious decomposition into
$3$ regular hexagons, each composed of $6$ equilateral triangles. The
$3$ edges of the complex associated to the centers of these three hexagons
have a common endpoint at $\infty$, and intersect $\C$ at three points
$p_1,p_2,p_3$. Six triangles come together at each of the points $p_i$, and
a horoball centered at each of these points intersects $6$ tetrahedra in a
hexagonal pattern. Glue opposite faces of these hexagons to produce $3$
cusps centered at $p_1,p_2,p_3$. The result of all this gluing produces
the manifold $M$. By construction, this particular choice of triangulation
produces $6$ hexagonal cusps, $3$ made up of $6$ equilateral triangles and
$3$ made up of $18$ equilateral triangles.

\begin{figure}[ht]
\scalebox{.35}{\includegraphics{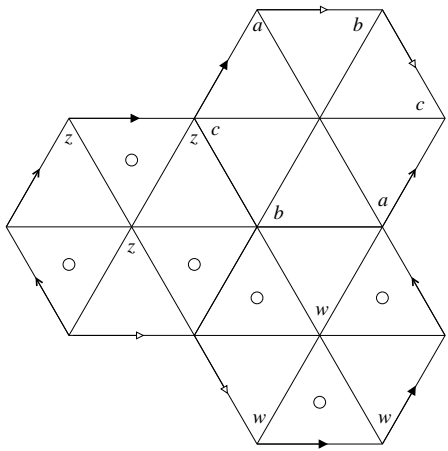}}
\caption{A horoball section of the cusp $d_1$ of Napoleon's $3$--manifold.
When Dehn surgery is performed on $c_1$, the triangles making up $d_1$
deform in the manner indicated. The triangles marked with a circle stay
equilateral under the deformation.}
\end{figure}

The components of $L$ fall into two
sets of $3$ links, depicted in the figure as the darker and the lighter
links, which we denote by $c_1,c_2,c_3$ and $d_1,d_2,d_3$
respectively. The group of symmetries of $M$ permutes the cusps by the
group $(\mathsf{S}_3 \times \mathsf{S}_3) \rtimes \Z/2\Z$ where the
conjugation action of the generator of $\Z/2\Z$ takes $(\sigma_i,\sigma_j)$
to $(\sigma_j,\sigma_i)$. If $G$ is the entire group of symmetries of
$M$, then there is a short exact sequence

$$ 0 \longrightarrow \Z/3\Z \longrightarrow G \longrightarrow
(\mathsf{S}_3 \times \mathsf{S}_3) \rtimes \Z/2\Z \longrightarrow 0 $$

\vskip 12pt

so the order of $G$ is $216$.

These links have the property that a surgery on $c_i$ keeps
$c_{i+1}$ and $c_{i+2}$ hexagonal, but distorts the structure at the
$d_j$, and vice versa. On the other hand, a surgery on {\em both} $c_i$ and
$c_{i+1}$ {\em does} distort the structure at $c_{i+2}$. In terms of the
picture already described, the configuration of $18$ triangles in
figure 4 decomposes into 3 hexagons of 6 triangles. These three hexagons
glue up in the obvious way to give hexagonal triangulations of the $c_i$.
Under surgery on $c_i$, the tetrahedra intersecting $c_i$ deform to
satisfy the new modified cusp equations. Under this deformation, the other
triangle types must deform to keep the other cusps complete. It can be
easily checked that the triangulations of the cusps $c_{i+1},c_{i+2}$ lift
to the symmetric tessellations of $\R^2$ depicted in figure 1, and therefore
the similarity types of these cusps stay hexagonal. However, under
surgery on {\em both} 
$c_i$ and $c_{i+1}$, the similarity types of triangles making
up cusp $c_{i+2}$ are not related in any immediately apparent way to
the picture in figure 1. The proof that $c_2$ and $c_3$ are 
isolated from $c_1$ 
is essentially just a calculation that under surgery on $c_1$ say,
the simplex parameters solving the relevant edge and cusp equations 
for the combinatorial triangulation vary as indicated in figure 4.
One may check experimentally that the similarity
type of $c_3$ is not constant under fillings on both
$c_1$ and $c_2$, using for example, Jeff Weeks' program
{\tt snappea}, available from \cite{jWsn}, for finding hyperbolic 
structures on $3$--manifolds, or Andrew Casson's provably accurate
program {\tt cusp}
(\cite{aCcu}).

\begin{defn}
Say that $3$ cusps are in {\em Brunnian isolation} when a surgery on one
of them leaves invariant the structure at the other $2$, but a surgery on
two of the cusps can change the structure of the third.
\end{defn}

With this definition, we observe that Napoleon's $3$--manifold has two sets
of cusps in Brunnian isolation.

One can see that there is an automorphism of $S^3$ of order $2$ fixing
two components of $L$ and permuting the other $4$ components in pairs.
The quotient by this automorphism is an orbifold $N$ which
has two regular cusps and two pillow cusps. We call the pillow cusps
$c_p, d_p$, and the regular cusps $c_r,d_r$ where $c_p$ is a quotient
of $c_1$, $c_r$ is covered by $c_2 \cup c_3$, and similarly for the $d_i$.

The cusp $c_p$ is first-order isolated but not isolated from $c_r$. 
Similarly, $d_p$ is first-order isolated but not isolated from $d_r$.
This can be easily observed by noting that $c_r$ {\em is} isolated from
$c_p$, and therefore it is first-order isolated from it (by the
corresponding properties in the cover). It follows
that $c_p$ is first-order isolated from $c_r$, by \cite{wNaR93}. To see
that $c_p$ is {\em not} isolated from $c_r$, it suffices to pass to the
cover and perform an equivariant surgery there. Alternatively, one can
easily check by hand using {\tt snappea} that the geometry of the 
cusp $c_p$ changes when one performs surgery on $c_r$.
\end{exa}

\begin{exa}
The $2$--cusped orbifold $A$ first described in \cite{wNaR93} and studied in
\cite{dC96} displays Napoleonic tendencies, where the version of
Napoleon's theorem we use now concerns right triangles. It is obtained by
$(2,0)$ surgery on the light cusps of the link complement portrayed in
figure 5. Coincidentally, it is $2$--fold covered by that very link
complement. Let $T$ be the right triangle with side lengths 
$\lbrace 1,1,\sqrt{2} \rbrace$, and let $S$ be the unit square in $\C$.
Pick a point $p \in S$ and construct four triangles $T_i$ all similar to
$T$ with one vertex of the diagonal at $p$ and the other vertex of the
diagonal at a vertex of $S$, such that the triangle is clockwise of the
diagonal, seen from $p$. Then the $8$ vertices of these triangles away
from $p$ are the vertices of an octagon which tiles the plane with
quotient space the square torus. This corresponds exactly to the
triangulation of a horoball section of a complete cusp in the orbifold $A$
after a deformation of the other cusp. One observes that the link 
complement in figure 3 is obtained from the link complement in figure 5
by drilling out two curves. Are the isolation phenomena associated with
the two links related?
\end{exa}

\begin{figure}[ht]
\scalebox{.4}{\includegraphics{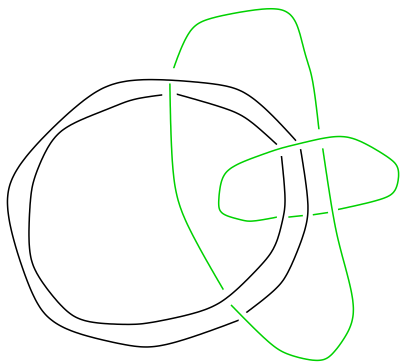}}
\caption{A link in $S^3$ whose complement displays Napoleonic tendencies}
\end{figure}

\end{document}